\documentclass[12pt]{article} 
\usepackage{amssymb}
\usepackage{amsmath}

\def\shadowbox{\hbox{\rule[-0.0ex]{0.1ex}{1.2ex}%
\hspace{-0.1ex}\rule[-0.0ex]{1.2ex}{0.1ex}%
\hspace{0.0ex}\rule[-0.0ex]{0.1ex}{1.2ex}\hspace{-1.3ex}%
\rule[1.15ex]{1.25ex}{0.1ex}\hspace{-0.0ex}\rule[-0.25ex]{0.3ex}{1.1ex}%
\hspace{-1.2ex}\rule[-0.25ex]{1.1ex}{0.25ex}}}
\def\qed{\ifmmode \hbox{\hfill\shadowbox}
     \else \hphantom{x}\hfill\shadowbox \fi}

\def\proof{\noindent{\bf Proof: }}

\newtheorem{theorem}{Theorem}[section]

\newtheorem{definition}[theorem]{Definition}

\newtheorem{corollary}[theorem]{Corollary}

\def\remark{{\medskip \noindent \bf Remark:\hspace{0.5em}}}

\def\Cst{{\mathbb C}}
\def\Est{{\mathbb E}}
\def\Nst{{\mathbb N}}

\def\Rst{{\mathbb R}}
\def\Tst{{\mathbb T}}
\def\Zst{{\mathbb Z}}
\def\eps{\varepsilon}
\def\theta{\vartheta}
\def\phi{\varphi}
\def\Lsp{{\boldsymbol L}}
\def\Fsp{{\cal F}}

\def\Usp{{\cal U}}
\def\Hsp{{\cal H}}
\def\Hspm{{\cal H}^m}
\def\Hspn{{\cal H}^n}

\def\Ltsp{{\Lsp^2}}

\def\LtR{{\Lsp^2(\Rst)}}

\def\grass{{\cal G}}
\def\gmn{{g_{m,n}}}

\def\glam{{g_{\lambda}}}
\def\glamp{{g_{\lambda'}}}
\def\ord{{\cal O}}

\def\trace{{\operatorname{trace}}}
\def\rank{{\operatorname{rank}}}
\def\argmin{{\operatorname{argmin}}}

\def\Ind{{\cal I}}

\def\Indo{{\cal I_0}}
\def\TT{{T^\ast}}
\def\tf{{\tilde{f}}}

\def\acos{{\operatorname{acos}}}
\def\vol{{\operatorname{vol}}}
\def\group{{\operatorname{group}}}
\def\family{{{\cal F}}}
\def\frame{{\{f_k\}_{k \in \Ind}}}
\def\frame{{\{f_k\}_{k \in \Ind}}}
\def\tightframe{{\{h_k\}_{k \in \Ind}}}
\def\finframe{{\{f_k\}_{k=1}^N}}
\def\infframe{{\{f_k\}_{k=1}^{\infty}}}
\def\gab{{\cal G}(g,\Lambda)}
\def\gabg{{\cal G}(\phi_{\sigma},\Lambda)}
\def\gauss{{\phi_{\sigma}}}
\def\gausslam{{\phi_{\lambda}}}

\def\gaussopt{{\phi_{\sigma}^o}}
\def\Lambdaopt{{\Lambda^{o}}}
\def\corr{{\cal M}}
\def\channel{{{\mathbf H}}}
\def\yerr{{\tilde y}}
\def\frec{{\tilde f}}
\def\erind{{\cal E}}
\def\erframe{{\{f_k\}_{k \in \remind}}}
\def\remind{{\cal R}}
\def\TE{{T_{\remind}}}
\def\TES{{T_{\remind}^{\ast}}}
\def\SI{{S^{-1}}}
\def\SQI{{S^{-\frac{1}{2}}}}

\begin{document}

\title{Grassmannian Frames with Applications to Coding and Communication}

\author{Thomas Strohmer\thanks{Department of Mathematics, University of
California, Davis, CA 95616-8633, USA; strohmer@math.ucdavis.edu. T.S. 
acknowledges support from NSF grant 9973373.} and
Robert Heath\thanks{Dept.\ of Electrical and Comp.\ Engineering 
The University of Texas at Austin, Austin, TX 78727, USA;
rheath@ece.utexas.edu}}
%and Arogyaswami Paulraj\thanks{Information Systems Laboratory, 
%Stanford University, Stanford, CA 94305-9510, USA.}} 

\date{}
\maketitle

%\vspace*{-6cm}
%\noindent
%{\small \tt Submitted to Appl.\ Comp.\ Harm.\ Anal.}

%\vspace*{6cm}

\begin{abstract}
For a given class $\Fsp$ of uniform frames of fixed redundancy we define
a Grassmannian frame as one that minimizes the 
maximal correlation $|\langle f_k,f_l \rangle|$ 
among all frames $\frame \in \Fsp$. We first analyze finite-dimensional 
Grassmannian frames. Using links to packings in Grassmannian spaces and 
antipodal spherical codes we derive bounds on the minimal achievable
correlation for Grassmannian frames. These bounds yield a simple condition 
under which Grassmannian frames coincide with uniform tight frames. We 
exploit connections to graph theory, equiangular line sets, and coding 
theory in order to derive explicit constructions of Grassmannian frames. 
Our findings extend recent results on uniform tight
frames. We then introduce infinite-dimensional Grassmannian frames
and analyze their connection to uniform tight frames for frames
which are generated by group-like unitary systems.
We derive an example of a Grassmannian Gabor frame by using connections
to sphere packing theory. Finally we discuss the application
of Grassmannian frames to wireless communication
and to multiple description coding.
\end{abstract}

%\noindent
%{\em AMS Subject Classification:} 

\noindent
{\em Key words:} Frame, Grassmannian spaces, spherical codes,
Gabor frame, multiple description coding, uniform tight frame,
conference matrix, equiangular line sets, unitary system.

\section{Introduction}
\label{intro}

Orthonormal bases are an ubiquitous and eminently powerful tool
that pervades all areas of mathematics. Sometimes however we find ourselves
in a situation where a representation of a function or an operator by an
overcomplete spanning system is preferable over the use of an orthonormal 
basis. One reason for this may be that an orthonormal basis with the
desired properties does not exist. A classical example occurs in Gabor
analysis, where the Balian-Low theorem tells us that orthonormal Gabor 
bases with good time-frequency localization cannot exist, while it is not 
difficult to find overcomplete Gabor systems with excellent time-frequency 
localization. Another important reason is the deliberate introduction of
redundancy for the purpose of error correction in coding theory. 

When dealing with overcomplete spanning systems one is naturally lead to
the concept of {\em frames}~\cite{Dau92}. Recall that a sequence of 
functions $\{f_k\}_{k \in \Ind}$ ($\Ind$ is a countable index set)
belonging to a
separable Hilbert space $\Hsp$ is said to be a {\em frame} for $\Hsp$
if there exist positive constants ({\em frame bounds)} $A$ and $B$ such that
\begin{equation}
A \|f\|_2^2 \le \sum_{k \in \Ind}|\langle f,f_k \rangle|^2 \le B \|f\|_2^2
\label{frame}
\end{equation}
for every $f \in \Hsp$. 

Even when there are good reasons to trade orthonormal bases for frames
we still want to preserve as many properties of orthonormal bases as possible.
There are many equivalent conditions to define an orthonormal basis 
$\{e_k\}_{k \in \Ind}$ for $\Hsp$, such as
\begin{equation}
f = \sum_{k \in \Ind}\langle f, e_k \rangle e_k, \quad 
\forall f \in \Hsp, \quad\text{and $\|e_k\|=1,\enspace \forall k\in \Ind$,}
\label{defonb1} 
\end{equation}
or
\begin{equation}
\text{$\{e_k\}_{k \in \Ind}$ is complete in $\Hsp$ and 
$\langle e_k,e_l \rangle = \delta_{k,l}$},
\label{defonb2}
\end{equation}
where $\delta_{k,l}$ denotes the Kronecker delta.

These two definitions suggest two ways to construct 
frames that are ``as close as possible'' to orthonormal bases. Focusing
on condition~\eqref{defonb1} we are naturally lead to uniform tight 
frames, which satisfy 
\begin{equation}
f = \frac{1}{A}\sum_{k \in \Ind}\langle f , f_k \rangle f_k, \enspace
\forall f \in \Hsp, \quad \text{and $\|f_k\|_2 = 1, \enspace 
\forall k \in \Ind$}, 
\label{utf}
\end{equation}
where $A$ is the lower frame bound. This class of frames has been frequently 
studied and is fairly well understood~\cite{Dau92,CK02,Han01,JS00,GKK01}.

As an alternative, as proposed in this paper, we focus on
condition~\eqref{defonb2}, which essentially states that the elements
of an orthonormal basis are perfectly uncorrelated. This suggests to search for
frames $\frame$ such that the maximal correlation 
$|\langle f_k, f_l \rangle|$ for all $k,l \in \Ind$ with $k \neq l$,
is as small as possible. This idea will lead us to so-called 
{\em Grassmannian frames}, which are characterized by the property that the 
frame elements have minimal cross-correlation among a given class of frames.
The name ``Grassmannian frames'' is motivated by the fact that
in finite dimensions Grassmannian frames coincide with 
optimal packings in certain Grassmannian spaces as we will see
in Section~\ref{s:finite}.

Recent literature on finite-dimensional frames~\cite{GKK01,CK02,EB02} 
indicates that the connection between finite frames and areas such as
spherical codes, algebraic geometry, graph theory, and sphere packings
is not well known in the ``frame community''. This has led to a number of
rediscoveries of classical constructions and duplicate results.
The concept of Grassmannian frames will allow us to make many of these
connections transparent.

\medskip

The paper is organized as follows. In the remainder of this
section we introduce some notation used throughout the paper.
In Section~\ref{s:finite} we focus on
finite Grassmannian frames. By utilizing a link to spherical codes
and algebraic geometry we derive lower bounds on the minimal achievable
correlation between frame elements depending on the redundancy of the
frame. We further show that optimal finite Grassmannian frames which
achieve this bound are also tight and certain uniform tight
frames are also Grassmannian frames. We discuss related concepts arising
in graph theory, algebraic geometry and coding theory and provide explicit 
constructions of finite Grassmannian frames.
In Section~\ref{s:grass} we extend the concept of Grassmannian frames to
infinite-dimensional Hilbert spaces and analyze the connection to 
uniform tight frames. We give an example of a Grassmannian
frame arising in Gabor analysis. Finally, in Section~\ref{s:appl} we discuss 
applications in wireless communication and coding theory.

\subsection{Notation} \label{ss:notation}

We introduce some notation and definitions used throughout the paper.
Let $\frame$ be a frame for a finite- or infinite-dimensional Hilbert space
$\Hsp$. Here $\Ind$ is an index set such as $\Zst, \Nst$ or $\{0,\dots,N-1\}$.
The {\em frame operator} $S$ associated with the frame $\frame$ is defined by
\begin{equation}
Sf = \sum_{k \in \Ind}\langle f,f_k \rangle f_k .
\label{frameop}
\end{equation}
$S$ is a positive definite, invertible operator that satisfies
$A I \le S \le B I$, where $I$ is the identity operator on $\Hsp$. The 
{\em frame analysis operator} $T: \Hsp \rightarrow \ell_2(\Ind)$ is given by
\begin{equation}
Tf = \{\langle f, f_k \rangle\}_{k \in \Ind}, 
\label{frameanalysis}
\end{equation}
and the {\em frame synthesis operator} is
\begin{equation}
T^{\ast}: \ell_2(\Ind) \rightarrow \Hsp: T \{c_{k}\}_{k \in \Ind}
= \sum_{k \in \Ind} c_k f_k.
\label{framesynthesis}
\end{equation}

Any $f \in \Hsp$ can be expressed as
\begin{equation}
f = \sum_{k \in \Ind} \langle f ,f_k \rangle h_k =
 \sum_{k \in \Ind} \langle f ,h_k \rangle f_k,
\label{frameexp}
\end{equation}
where $\{h_k\}_{k \in \Ind}$ is the {\em canonical dual frame}
given by $h_k = S^{-1}f_k$.
If $A=B$ the frame is called {\em tight}, in which case $S = AI$ and
$h_k = \frac{1}{A}f_k$. The tight frame canonically associated to
$\{f_k\}_{k \in I}$ is $S^{-\frac{1}{2}} f_k$. 

If $\|f_k\| = 1$ for all $k$ then $\frame$ is called a 
{\em uniform frame}. Here $\|.\|$ denotes the $\ell_2$-norm of a vector
in the corresponding finite- or infinite-dimensional Hilbert space.
{\em Uniform tight frames} have many nice properties which make them an 
important tool in theory~\cite{RT98,Hei02} and in a variety of 
applications~\cite{GKK01,SB01,Str01,FS98}. Observe that if $\frame$
is a uniform frame, then $\{\SQI f_k\}_{k \in \Ind}$ is a tight frame,
but in general no longer uniform!

We call a uniform frame $\frame$ {\em equiangular} if
\begin{equation}
\label{defequi}
|\langle f_k,f_l \rangle| = c \quad 
\text{for all $k,l$ with $k \neq l,$}
\end{equation}
for some constant $c\ge 0$.
Obviously any orthonormal basis is equiangular.

\section{Finite Grassmannian frames, spherical codes, and equiangular lines} 
\label{s:finite}

In this section we concentrate of frames $\finframe$ for $\Est^m$ where
$\Est = \Rst$ 
or $\Cst$. As mentioned in the introduction we want to construct frames
$\finframe$ such that the maximal correlation 
$|\langle f_k, f_l \rangle|$ for all $k,l \in \Ind$ with $k \neq l$,
is as small as possible. If we do not impose any other conditions on the
frame we can set $N=m$ and take $\finframe$ to be an 
orthonormal basis. But if we want to go beyond this trivial case and
assume that the frame is indeed overcomplete then the correlation
$|\langle f_k, f_l \rangle|$ will strongly depend on the redundancy
of the frame, which can be thought of as a 
``measure of overcompleteness''. Clearly, the smaller the redundancy the 
smaller we expect $|\langle f_k, f_l \rangle|$ to be.
In $\Est^m$ the redundancy $\rho$ of a frame $\finframe$ is defined by
$\rho = \frac{N}{m}$. 

\begin{definition}
For a given uniform frame $\finframe$ in $\Est^m$ we define 
%the {\em frame correlation} $\ave(\finframe)$ by
%\begin{equation}
%\ave(\finframe)=\sum_{k,l=1}^{N}|\langle f_k,f_l \rangle|^2,
%\label{ave}
%\end{equation}
the {\em maximal frame correlation} $\corr(\finframe)$ by
\begin{equation}
\corr(\finframe)=\underset{k,l, k\neq l}{\max}\{|\langle f_k,f_l \rangle|\}.
\label{corr}
\end{equation}
\end{definition}
The restriction to uniform frames in the definition above
is just for convenience, alternatively we could consider general frames
and normalize the inner product in~\eqref{corr} by the 
norm of the frame elements. Hence without loss of generality we can assume 
throughout this section that all frames are uniform.

\begin{definition}
\label{def:fin}
A sequence of vectors $\{u_k\}_{k=1}^{N}$
in $\Est^m$ is called a {\em Grassmannian frame} if it is the solution to
\begin{equation}
\min \big\{ \corr(\finframe) \big\},
%\min \{ \underset{k,l, k \neq l}{\max} \{ |\langle f_k, f_l \rangle|\},
\label{fingrass}
\end{equation}
where the minimum is taken over all uniform frames $\{f_k\}_{k=1}^{N}$ in
$\Est^m$.
\end{definition}

In other words a Grassmannian frame minimizes the maximal correlation
between frame elements among all uniform frames which have the same
redundancy. Obviously the minimum in~\eqref{fingrass} depends only on the 
parameters $N$ and $m$.

Two problems arise naturally when studying finite Grassmannian frames:  \\
{\bf Problem 1:} Can we derive bounds on $\corr(\finframe)$ for 
given $N$ and $m$? \\
{\bf Problem 2:} How can we construct Grassmannian frames? 

The following theorem provides an exhaustive answer to problem~1.
%It is the main theorem for finite Grassmannian frames.
The theorem is new in frame theory but actually it only unifies and
summarizes results from various quite different research areas. 
\begin{theorem}
\label{th:bound}
Let $\finframe$ be a frame for $\Est^m$. Then
\begin{equation}
\corr(\finframe) \ge \sqrt{\frac{N-m}{m(N-1)}}.
\label{bound1}
\end{equation}
Equality holds in~\eqref{bound1} if and only if
\begin{equation}
\label{equi}
\text{$\finframe$ is an equiangular tight frame.}
\end{equation}
Furthermore, 
\begin{equation}
\label{numbound1}
\text{if $\Est=\Rst$ equality in~\eqref{bound1} can only hold
if}\,\, N \le \frac{m(m+1)}{2}, 
\end{equation}
\begin{equation}
\label{numbound2}
\text{if $\Est=\Cst$ then equality in~\eqref{bound1} can only hold if 
$N \le m^2$.}
\end{equation}
\end{theorem}

\begin{proof}
A proof of the bound~\eqref{bound1} can be found in~\cite{Wel74,Ros97}. It 
also follows from Lemma~6.1 in~\cite{LS66}. One way to derive~\eqref{bound1}
is to consider the non-zero eigenvalues $\lambda_1,\dots,\lambda_m$
of the Gram matrix $R = \{\langle f_k,f_l \rangle\}_{k,l =1}^N$.
These eigenvalues satisfy $\sum_{k=1}^{m} \lambda_k = N$ and also
\begin{equation}
\label{lam2}
\sum_{k=1}^{m} \lambda_k^2 =
\sum_{k=1}^{N}\sum_{l=1}^{N}|\langle f_k,f_l \rangle|^2 \ge \frac{N^2}{m},
\end{equation}
see~\cite{Ros97,LS66}.
The bound follows now by taking the maximum over 
all $|\langle f_k,f_l \rangle|$ in~\eqref{lam2} and 
observing that there are $N(N-1)/2$ different pairs 
$\langle f_k,f_l \rangle$ for $k,l \in \Ind$ with $k\neq l$.

Equality in~\eqref{bound1} implies 
$\lambda_1 = \dots = \lambda_m = \frac{N}{m}$,
which in turn implies tightness of the frame, and
also $|\langle f_k, f_l \rangle|^2 = \frac{N-m}{m(N-1)}$ for
all $k,l$ with $k \neq l$ which yields the equiangularity
(cf.~also \cite{Ros97,CHS96}). Finally the bounds on $N$ 
in~\eqref{numbound1},~\eqref{numbound2} follow from the 
bounds in Table II of~\cite{DGS75}.
\end{proof}

We call uniform frames that achieve the bound~\eqref{bound1}
{\em optimal Grassmannian frames}. The following corollary will be
instrumental in the construction of a variety of optimal 
Grassmannian frames.
\begin{corollary}
\label{corgram}
Let $m,N \in \Nst$ with $N \ge m$. Assume $R$ is a hermitian $N\times N$ 
matrix with entries $R_{k,k}=1$ and 
\begin{equation}
\label{Grammatrix}
R_{k,l} = 
\begin{cases}
\pm \sqrt{\frac{N-m}{m(N-1)}}, & \text{if $\Est = \Rst$,} \\
\pm i\sqrt{\frac{N-m}{m(N-1)}}, & \text{if $\Est = \Cst$,}
\end{cases}
\end{equation}
for $k,l=1,\dots,N; k\neq l$.
If the eigenvalues $\lambda_1,\dots,\lambda_N$ of $R$ are such that
$\lambda_1 = \dots = \lambda_m = \frac{N}{m}$ and
$\lambda_{m+1} = \dots = \lambda_N = 0$, then there exists
a frame $\finframe$ in $\Est^m$ that achieves the bound~\eqref{bound1}.
\end{corollary}
\begin{proof}
Since $R$ is hermitian it has a spectral factorization of the form
$R = W \Lambda W^{\ast}$,
where the columns of $W$ are the eigenvectors and the diagonal
matrix $\Lambda$ contains the eigenvalues of $R$. Without loss of
generality we can assume that the non-zero eigenvalues of $R$ are contained
in the first $m$ diagonal entries of $\Lambda$. 
Set $f_k := \sqrt{\frac{N}{m}} \{W_{k,l}\}_{l=1}^{m}$ for $k=1,\dots,N$. 
By construction we have $\langle f_k,f_l \rangle = R_{l,k}$, hence $\finframe$
is equiangular. Obviously $\finframe$ is tight, since all non-zero
eigenvalues of $R$ are identical. Hence by Theorem~\ref{th:bound} $\finframe$
achieves the bound~\eqref{bound1}. 
\end{proof}

On the first glance Corollary~\ref{corgram} does not seem to make the
problem of constructing optimal Grassmannian frames much easier.
However by using a link to graph theory and spherical designs we 
will be able to derive many explicit constructions of matrices having
the properties outlined in Corollary~\ref{corgram}.

While the concept of Grassmannian frames is new in frame theory there are 
a number of related concepts in other areas of mathematics.
Thus it is time to take a quick journey through these areas which
will take us from Grassmannian spaces to spherical designs to coding theory.

\medskip
\noindent
{\bf Packings in Grassmannian spaces:}

The {\em Grassmannian space} $\grass(m,n)$ is the set of all
$n$-dimensional subspaces of the space $\Rst^m$ 
(usually the Grassmannian space is defined for $\Rst$ only, although many
problems can be analogously formulated for the complex space).
$\grass(m,n)$ is a homogeneous space isomorphic to 
$\ord(m)/(\ord(n) \times \ord(m-n))$, it forms a compact Riemannian 
manifold of dimension $n(m-n)$. 

The {\em Grassmannian packing problem} is the problem of finding
the best packing of $N$ $n$-dimensional subspaces in $\Est^m$, such that
the angle between any two of these subspaces becomes as large as possible
\cite{CHS96,CHR99}.
In other words, we want to find N points in $\grass(m,n)$ so that the minimal
distance between any two of them is as large as possible. 
For our purposes we can concentrate on the case $n=1$. Thus the subspaces
are (real or complex) lines through the origin in $\Est^m$ and the goal is to 
arrange $N$ lines such that the angle between any two of the lines becomes
as large as possible. Since maximizing the angle between lines is
equivalent to minimizing the modulus of the inner product of the
unit vectors generating these lines, it is obvious that finding
optimal packings in $\grass(m,1)$ is equivalent to finding finite Grassmannian
frames (which also motivated the name for this class of frames).

By embedding the Grassmannian space $\grass(m,n)$ 
into a sphere of radius $\sqrt{n(m-n)/m}$ in $\Rst^d$ with $d=(m+1)m/2-1$, 
Conway, Hardin, and Sloane are able to apply bounds from spherical codes 
due to Rankin~\cite{Ran55} to derive bounds on the maximal angle between 
$N$ subspaces in $\grass(n,m)$, see the very inspiring paper~\cite{CHS96}.
For the case $n=1$ the bound coincides of course with~\eqref{bound1}.

\medskip
\noindent
{\bf Spherical codes:}

A {\em spherical code} ${\cal S}(m,N,s)$ is a set of $N$ points (code words) 
on the $m$-dimensional unit sphere $\Omega_m$, such that the inner product
between any two code words is smaller than $s$, cf.~\cite{CS93b}.
By placing the points on the sphere as far as possible from each other
one attempts to minimize the risk of decoding errors.
{\em Antipodal spherical codes} are spherical codes
which contain with each code word $w$ also the code word $-w$. 
Clearly, the construction of antipodal spherical codes whose $N$ points
are as from each other as possible is closely related to constructing
Grassmannian frames.

In coding theory the inequality at the right-hand side of~\eqref{lam2} is 
known as {\em Welch bound}, cf.~\cite{Wel74}. Therefore uniform tight frames
are known as {\em Welch bound equality} (WBE) sequences in coding
theory\footnote{The authors of~\cite{VAT99} incorrectly call WBE sequences 
tight frames.}. WBE sequences have gained new popularity in connection with
the construction of spreading sequences for Code-Division Multiple-Access
(CDMA) systems~\cite{VAT99,HS02,Sar99}. WBE sequences that
meet~\eqref{bound1} with equality are called maximum WBE (MWBE)
sequences~\cite{Wel74,Sar99}. While Welch (among other authors)
derived the bound~\eqref{bound1} he did not give an explicit construction
of MWBE sequences. 

\medskip
\noindent
{\bf Spherical designs:}

A {\em spherical $t$-design}\footnote{A {\em spherical} $t$-design should
not be confused with an ``ordinary'' $t$-design.} is a finite subset $X$ of
the unit sphere $\Omega_m$ in $\Rst^m$, such that
\begin{equation}
|X|^{-1} \sum_{x \in X} h(x) = \int \limits_{\Omega_m} h(x) dw(x),
\label{sphericaldesign}
\end{equation}
for all homogeneous polynomials $h \in \text{Hom}_t(\Rst^m)$ of total
degree $t$ in $m$ variables, see e.g.~\cite{Sei01}. A spherical design
measures certain regularity properties of sets $X$ on the unit sphere
$\Omega_m$. Another way to define a spherical $t$-design is by requiring
that, for $k=0,\dots,t$ the $k$-th moments of $X$ are constant with
respect to orthogonal transformations of $\Rst^m$. Here are a few
characterizations of spherical $t$-designs that make the connection
to the aforementioned areas transparent. For details about the following
examples we refer to~\cite{DGS77}. Let the cardinality of $X$ be $N$.
$X$ is a spherical 1-design if and only if the Gram matrix $R(X)$
of the vectors of $X$ has vanishing row sums. $X$ is a spherical
2-design if it is a spherical 1-design and the Gram matrix $R(X)$ has
only two different eigenvalues, namely $N/m$ with multiplicity $m$, and
0 with multiplicity $N-m$. An antipodal spherical code on $\Omega_m$ is
a 3-design if and only if the Gram matrix of the corresponding set of
vectors has two eigenvalues.

\medskip
\noindent
{\bf Equiangular line sets and equilateral point sets:}

In~\cite{LS73a,DGS75} Seidel et al.\ consider sets of lines
in $\Rst^m$ and in $\Cst^m$ having a prescribed number of angles.
They derive upper bounds on the number of lines in the case of
one, two, and three prescribed angles (in the latter case, one of the
angles is assumed to be zero). Most interesting are those line sets that
actually meet the upper bound. 
In~\cite{LS66} van~Lint and Seidel consider a similar problem
in elliptic geometry. Since the unit sphere in $\Rst^m$ serves as model for 
the $m-1$-dimensional elliptic space ${\cal \Est}^{m-1}$ where any elliptic 
point is represented by a pair of antipodal points in $\Rst^m$, the 
construction of equilateral point sets in elliptic geometry is of course 
equivalent to the construction of equiangular lines sets in
Euclidean geometry. Recall that optimal Grassmannian frames
are equiangular, hence the search for equiangular line sets is closely
related to the search for optimal Grassmannian frames.

\medskip
\noindent
{\bf Characterization of strongly regular graphs:} 

Graphs with a lot of structure and symmetry play a central role in graph 
theory. Different kinds of matrices are used to represent a graph, such as 
the Laplace matrix or adjacency matrices~\cite{BCN89}. What 
structural properties can be derived from the eigenvalues depends on the 
specific matrix that is used. The Seidel adjacency matrix $A$ of
a graph $\Gamma$ is given by
\begin{equation}
A_{xy} = 
\begin{cases}
-1 & \text{if the vertices $x, y \in \Gamma$ are adjacent,} \\
 1 & \text{if the vertices $x, y \in \Gamma$ are non-adjacent,} \\
 0 & \text{if $x=y$.}
\end{cases}
\label{seidelmat}
\end{equation}

If $A$ has only very few different eigenvalues then the graph is (strongly)
regular, cf.~\cite{BCN89}. The connection to Grassmannian frames $\finframe$
that achieve the bound~\eqref{bound1} is as follows. Assume that the
associated Gram matrix $R = \{\langle f_k, f_l \rangle\}_{k,l=1}^N$ 
has entries $\pm \alpha$ and 1 at the diagonal. Then
\begin{equation}
A = \frac{1}{\alpha}(R-I)
\label{graphmat}
\end{equation}
is the adjacency matrix of a regular two-graph~\cite{Sei76}.
We will make use of this relation in the next section.

\subsection{Construction of optimal Grassmannian frames} \label{ss:}

In this section we give some explicit constructions for infinite
families of optimal finite Grassmannian frames. Note that
optimal Grassmannian frames do not exist for all choices of $m$
and $N$ (assuming of course that $N$ does not exceed $(m+1)m/2$ or $m^2$,
respectively). For instance there are no 5 vectors in $\Rst^3$ with maximal
correlation $\frac{1}{\sqrt{6}}$. In fact, although the 5 vectors in $\Rst^3$ 
that minimize~\eqref{fingrass} are equiangular, the maximal inner product is
$\frac{1}{\sqrt{5}}$ (but not $\frac{1}{\sqrt{6}}$),
see~\cite{CHS96}. On the other hand the 7 vectors in $\Rst^3$ that 
minimize~\eqref{fingrass} yield a uniform tight frame, but not an equiangular
one (which should not come as a surprise since the choice $N=7$ exceeds
the bound $N \le m(m+1)/2$. Note that for $\Cst^3$ we can indeed
construct 7 lines that achieve the bound~\eqref{bound1}.
We refer to~\cite{CHS96} for details about some of these and other examples. 

\begin{corollary}
\label{utfgrass}
Let $\Est = \Rst$ or $\Cst$ and $N=m+1$. $\finframe$ is an optimal 
Grassmannian frame for $\Est^m$ if and only if it is a uniform tight frame.
\end{corollary}
\begin{proof}
An optimal Grassmannian frame $\finframe$ with $N=m+1$ can be
easily constructed by taking the vectors to be the vertices of a regular 
simplex in $\Est^m$, cf~\cite{CHS96}. Thus by Theorem~\ref{th:bound} 
$\finframe$ is
a uniform tight frame. On the other hand it was shown in~\cite{GKK01}
that all uniform tight frames with $N=m+1$ are equivalent. Since
this equivalence relation preserves inner products it follows that
any uniform tight frame $\finframe$ with $N=m+1$ achieves the
bound~\eqref{bound1}.
\end{proof}

A uniform tight frame $\finframe$ with $N=m+1$ also provides a
spherical 1-design, which can be seen as follows. When $N=m+1$ we can 
always multiply the elements of $\finframe$ by $\pm 1$ such that
the Gram matrix $R$ has 1 as its main diagonal entries and $-\frac{1}{m}$
else. Hence the row sums of $R$ vanish and therefore $\finframe$
constitutes a spherical 1-design. It is obvious that the Seidel adjacency 
matrix $A$ of a graph which is constructed from a regular simplex
has $A_{kk}=0$ and $A_{kl}=-1$ for $k \neq l$, which illustrates nicely the 
relationship between $A$ and $R$ as stated in~\eqref{seidelmat}.

\if 0
According to~\cite{CHS96} the bound~\eqref{bound1} is essentially 
Rankin's bound on spherical codes~\cite{Ran55}, also known as 
{\em Simplex bound}.
It has been rediscovered (sometimes in a somewhat different form)
by researchers in areas ranging from coding
theory~\cite{Wel74}, to graph theory~\cite{DGS75} and algebraic 
geometry~\cite{Ros97}. A fascinating contribution
is due to Conway, Hardin, and Sloane~\cite{CHS96} who derive a 
bound for the maximal angle between $n$-dimensional subspaces in $\Rst^m$ 
(which reduces to the present bound if $n=1$)
by an elegant and very natural embedding approach.

The fact that Grassmannian frames which achieve the bound~\eqref{bound1}
are automatically tight is very useful in applications.

\noindent
{\bf Remark:} We can fix a value $\alpha$ and ask
how large can $N$ be such that $\rho(N,m) = \alpha$ or 
$\rho(N,m) \le \alpha$ still holds.
In other words what is the maximum cardinality $N_{\max}$
of a spherical code with mutual correlation $\alpha$?
From Theorem~\ref{th:bound} we get 
\begin{equation}
N_{\max} \le \frac{m(1-\alpha)}{1-m\alpha},
\label{maxN}
\end{equation}
as shown e.g.\ in~\cite{DGS75}. More generally we may ask e.g.\ what is 
$N_{\max}$ such that $|\langle u_k,u_l \rangle| \in \{\alpha,\beta\}$
for all $k,l$ with $k \neq l$. We refer to the work of Seidel et al.\
for a detailed investigations of this and related 
problems~\cite{DGS75,DGS77}. See also~\cite{HS02} 
in connection with out-of-cell interference for CDMA.
\fi

\bigskip
The following construction has been proposed in~\cite{LS73,CHS96}.
An $n \times n$ {\em conference matrix} $C$
has zeros along its main diagonal and $\pm 1$ as its other entries,
and satisfies $C C^T = (n-1)I_n$, see~\cite{LS73}.
Conference matrices play an important role in graph theory~\cite{Sei76}.
If $C_{2m}$ is a symmetric conference matrix, then there exist
exist $2m$ vectors in $\Rst^m$ such that the bound~\eqref{bound1} holds
with $\rho(2m,m) = 1/\sqrt{2m-1}$. 
If $C_{2m}$ is a skew-symmetric conference matrix
(i.e., $C=-C^T$), then there exist exist $2m$ vectors in $\Cst^m$ such 
that the bound~\eqref{bound1} holds with $\rho(2m,m) = 1/\sqrt{2m-1}$,
see Example~5.8 in~\cite{DGS75}. 
The link between the existence of a (real or complex) optimal
Grassmannian frame and the existence of a corresponding
conference matrix $C_{2m}$ can be easily as seen as follows.
Assume that $\finframe$ achieves~\eqref{bound1} and denote 
$\alpha:=1/\sqrt{2m-1}$. We first consider the case $\Est=\Rst$. Clearly
the entries of the $2m \times 2m$ Gram matrix 
$R = \{\langle f_k,f_l \rangle\}_{k,l=1}^{N}$ are 
$R_{k,l}=\pm \alpha$ for $k\neq l$ and $R_{k,k}=1$. Hence
\begin{equation}
\label{gram2conf1}
C:= \frac{1}{\alpha}(R-I)
\end{equation}
is a symmetric conference matrix.
For $\Est=\Cst$ we assume that $R_{k,l}=\pm i\alpha$
for $k\neq l$ and $R_{k,k}=1$. Then
\begin{equation}
\label{gram2conf2}
C:= \frac{1}{i\alpha}(R-I)
\end{equation}
is a skew-symmetric conference matrix.

The derivations above lead to the following
\begin{corollary}
\label{confgrass}
(a) Let $N=2m$, with $N=p^{\alpha}+1$ where $p$ is a prime number and 
$\alpha \in \Nst$. Then there exists an optimal Grassmannian frame
in $\Rst^m$ which can be constructed explicitly. \\
(b) Let $N=2m$, with $m=2^\alpha$ with $\alpha \in \Nst$. Then there exists 
an optimal Grassmannian frame in $\Cst^m$ which can be constructed explicitly. 
\end{corollary}
\proof
Paley has shown that if $N = p^\alpha +1$ with
$p$ and $\alpha$ as stated above, then there exists a symmetric
$N \times N$ conference matrix, moreover this matrix can be constructed
explicitly, see~\cite{Pal33,GS67}. For the case $N=2m=2^{\alpha+1}$ a
skew-symmetric conference matrix can be constructed by the following 
recursion: Initialize
\begin{equation}
C_2 = \begin{bmatrix}
0 & -1 \\
1 & 0
\end{bmatrix},
\label{conf2}
\end{equation}
and compute recursively
\begin{equation}
C_{2m} = \begin{bmatrix}
C_m & C_m - I_m \\
C_m + I_m & -C_m
\end{bmatrix},
\label{conf2m}
\end{equation}
then it is easy to see that $C_{2m}$ is a skew-symmetric conference matrix.

An application of Corollary~\ref{corgram} to both, the symmetric and the
skew-symmetric conference matrix respectively, completes the proof.
\qed

Hence for instance there exist 50 equiangular lines in $\Rst^{25}$ with angle
$\acos (1/\sqrt{49})$ and 128 equiangular lines in $\Cst^{64}$ with angle 
$\acos (1/\sqrt{127})$. The construction in~\eqref{conf2},~\eqref{conf2m} is 
reminiscent of the construction of Hadamard matrices. Indeed, $C_m + I_m$ is a 
skew-symmetric Hadamard matrix.

\subsubsection{Nearly optimal Grassmannian frames} \label{sss:nearopt}
\bigskip

Theorem~\ref{th:bound} gives an upper bound on the cardinality of
optimal Grassmannian frames. If the redundancy of a frame is too large
then it cannot achieve equality in~\eqref{bound1}. But it is possible
to design Grassmannian frames whose cardinality slightly exceeds
the bounds in Theorem~\ref{th:bound}, while their maximal correlation
is close to the optimal value. For instance there exist frames
$\finframe$ in $\Cst^m$ where $N=m^2+1$, with maximal correlation
$\corr = \frac{1}{\sqrt{m}}$. In fact, these nearly optimal Grassmannian
frames are unions of orthonormal bases, and the modulus of the
inner products between frame elements takes on only the values $0$
and $\frac{1}{\sqrt{m}}$. We refer to~\cite{CCK97,Lev82} for details
about these amazing constructions, which find an important
application in the design of spreading sequences for CDMA~\cite{HS02}.

\medskip
\noindent
{\bf Example:}
Here is an example of a discrete finite Gabor frame that is a nearly
optimal Grassmannian frame in $\Cst^m$ (see~\cite{FS98} for details
about Gabor frames).
Let $m$ be a prime number $\ge 5$ and set $g(n) = e^{2\pi i n^3/m}$
for $n=0,\dots,m-1$. Then the frame $\{g_{k,l}\}_{k,l=0}^{m-1}$, where
\begin{equation}
g_{k,l}(n) = g(n-k) e^{2\pi i ln/m}, \quad k,l=0,\dots,m-1,
\label{fingabor}
\end{equation}
satisfies $|\langle g_{k,l}, g_{k',l'} \rangle| \in \{0,1/\sqrt{m}\}$
for all $g_{k,l} \neq g_{k',l'}$, which follows from Theorem~2
in~\cite{All80}. Hence ${\cal M}(\{g_{k,l}\}_{k,l=0}^{m-1})=1/\sqrt{m}$ while
\eqref{bound1} yields $1/\sqrt{m+1}$ as theoretically optimal value.
Note that we can add the standard orthonormal basis to the frame 
$\{g_{k,l}\}_{k,l=0}^{m-1}$ without changing the maximal frame correlation
$1/\sqrt{m}$.

\section{Infinite-dimensional Grassmannian frames} \label{s:grass}

In this section we extend the concept of Grassmannian frames to
frames $\infframe$ in separable infinite-dimensional Hilbert spaces. 
As already pointed out in Section~\ref{s:finite} the maximal
correlation $|\langle f_k,f_l \rangle|$
of the frame elements will depend crucially on the redundancy of the frame.
While it is clear how to define redundancy for finite frames,
it is less obvious for infinite dimensional frames.

The following appealing definition is due to 
Radu Balan and Zeph Landau~\cite{BL02}.
\begin{definition}
\label{def:red}
Let $\{f_k\}_{k=1}^{\infty}$ be a frame for $\Hsp$. The redundancy $\rho$
of the frame $\{f_k\}_{k=1}^{\infty}$ is defined as
\begin{equation}
\rho := \Big(\lim_{K \rightarrow \infty} \frac{1}{K}
\sum_{k=1}^{K}
\langle f_k, S^{-1}f_k \rangle\Big)^{-1},
\label{redundancy}
\end{equation}
provided that the limit exists.
\end{definition}
Using the concept of ultrafilters Balan and Landau have derived a more
general definition of redundancy of frames, which coincides of course with the
definition above whenever the limit in~\eqref{redundancy} exists~\cite{BL02}. 
In this paper we will restrict ourselves to the definition of redundancy 
as stated in~\eqref{redundancy} since it is sufficiently general for
our purposes.

\remark
We briefly verify that the definition of frame redundancy by Balan and Landau 
coincides with our usual understanding of redundancy in some
important special cases:\\
(i) Let $\{f_k\}_{k =1}^N$ be a finite frame for an $m$-dimensional Hilbert 
space $\Hspm$. Let $P: \Hspn \rightarrow \Hspm$ denote the associated
projection matrix with entries $P_{k,l} = \langle f_k,S^{-1}f_l \rangle$. 
We compute 
\begin{equation}
\label{red1}
\rho=\Big(\frac{1}{N}\sum_{k=1}^{N}\langle f_k, S^{-1}f_k \rangle\Big)^{-1}= 
\frac{N}{\trace (P)} = \frac{N}{\rank (P)}= \frac{N}{m},
\end{equation}
which coincides with the usual definition of redundancy in finite dimensions.

(ii) Let $\{\gmn\}_{m,n \in \Zst}$, where $\gmn(x) = g(x-ma) e^{-2\pi i nbx}$
be a Gabor frame for $\LtR$ with time- and frequency-shift parameters 
$a,b >0$.  We have from~\cite{Jan95} that
\begin{equation}
\label{red2}
\langle \gmn,S^{-1}\gmn \rangle = \langle g, S^{-1}g \rangle = ab,
\qquad \text{for all $m,n \in \Zst$}, 
\end{equation}
hence $\rho = 1/(ab)$ as expected. 

(iii) Assume $\frame$ is a uniform tight frame. Then 
$\langle f_k, S^{-1} f_k \rangle = 1/A$ and therefore
$\rho = A$, which agrees with the intuitive expectation that for
uniform tight frames the frame bound measures the redundancy
of the frame~\cite{Dau92}.

We need two more definitions before we can introduce the concept of
Grassmannian frames in infinite dimensions.

\begin{definition}[\cite{DL98}]
A {\em unitary system} $\Usp$ is a countable set of unitary operators
containing the identity operator and acting on a separable Hilbert space
$\Hsp$.
\end{definition}

\begin{definition}
\label{def:grassframe}
Let $\Usp$ be a unitary system and $\Phi$ be a class of 
functions with $\|\phi\|_2 =1$ for $\phi \in \Phi$. We denote by 
$\Fsp(\Hsp,\Usp,\Phi)$ the family of frames $\{\phi_k\}_{k \in \Ind}$ 
for $\Hsp$ of fixed redundancy $\rho$, such that
\begin{equation}
\label{grassdef1}
\phi_k = U_k \phi_0, \qquad \phi_0 \in \Phi, \enspace
U_k \in \Usp, \enspace k \in \Ind.
\end{equation}
We say that $\{f_k\}_{k \in \Ind}$ is a Grassmannian frame with respect to
$\Fsp(\Hsp,\Usp,\Phi)$ if it is the solution of~\footnote{For a 
frame $\{\phi_k\}_{k \in \Ind}$ there always exists
$\underset{k \neq l}{\max}\{|\langle \phi_k, \phi_l \rangle|\}$, 
otherwise the upper frame bound could not be finite.}
\begin{equation}
\label{grassdef2}
%\underset{\{\phi_k\}_{k \in \Ind} \in \Fsp(\Hsp,\Usp,\Phi)}{\argmin}
%\Big(\underset{k,l \in \Ind; k \neq l}{\max}
\underset{\{\phi_k\}_{k \in \Ind} \in \Fsp(\Hsp,\Usp,\Phi)}{\argmin}
\Big(\underset{k,l \in \Ind; k \neq l}{\max}
\big\{|\langle \phi_k, \phi_l \rangle| \big\}\Big)
\end{equation}
for given $\rho$.
\end{definition}

In the definition above we have deliberately chosen $\Phi$
such that it does not necessarily have to coincide with all
functions in $\Ltsp(\Hsp)$. The reason is that in many applications one 
is interested in designing frames using only a specific class of functions.

In finite dimensions we derived conditions under which Grassmannian
frames are also uniform tight frames. Such a nice and simple relationship 
does not exist in infinite dimensions. However in many cases it is possible
to construct a uniform tight frame whose maximal frame correlation 
is close to that of a Grassmannian frame as we will see in the
next theorem.

The following definition is due to Gabardo and Han~\cite{GH01}.
\begin{definition}
Let $\Tst$ denote the circle group. A unitary system $\Usp$ is called
{\em group-like} if
\begin{equation}
\group(\Usp) \subset \Tst \Usp : = \{tU \, : \, t \in \Tst, U \in \Usp\},
\label{grouplike}
\end{equation}
and if different $U, V \in \Usp$ are always linearly independent,
where $\group(\Usp)$ denotes the group generated by $\Usp$.
\end{definition}

\begin{theorem}
\label{th:tight}
Let $\Fsp(\Hsp,\Usp,\Phi)$ be given, where $\Usp$ is a group-like 
unitary system. For given redundancy $\rho$ assume that $\frame$ is 
a Grassmannian frame for $\Fsp(\Hsp,\Usp,\Phi)$ with frame bounds $A, B$. 
Then there exists a uniform tight frame $\tightframe$ with 
$h_k = U_k h_0, U_k\in \Usp$, such that
\begin{equation}
\underset{k,l \in \Ind; k \neq l}{\max} |\langle h_k,h_l \rangle|
\le \underset{k,l \in \Ind; k \neq l}{\max}|\langle f_k,f_l \rangle| +
     2 \max\Big\{|1-\sqrt{\frac{\rho}{A}}|,|1-\sqrt{\frac{\rho}{B}}|\Big\}.
\label{tightgrass}
\end{equation}
\end{theorem}
\begin{proof}
Let $S$ be the frame operator associated with the Grassmannian
frame $\frame$. We define the tight frame $\{h_k\}_{k \in \Ind}$ via
$h_k : = \sqrt{\rho}\SQI f_k$.
Since $\Usp$ is a group-like unitary system it follows
from~\eqref{grassdef1} above and Theorem~1.2 in~\cite{Han01} that 
\begin{equation}
\label{inf0}
\langle f_k, \SI f_k \rangle = \langle U_k f_0, \SI U_k f_0 \rangle
= \langle U_k f_0, U_k \SI f_0 \rangle = \langle f_0, \SI f_0 \rangle,
\end{equation}
and
\begin{equation}
\label{inf1}
h_k =\sqrt{\rho} \SQI f_k =\sqrt{\rho} \SQI U_k f_0 =\sqrt{\rho} U_k \SQI f_0.
\end{equation}
Using Definition~\ref{def:red} and~\eqref{inf0}, we get
$\langle f_k, \SI f_k \rangle = \frac{1}{\rho}$ and therefore
\begin{equation}
\label{inf2}
\|h_k\|^2 = \rho\langle \SQI f_k,\SQI f_k \rangle
=\rho \langle f_k,\SI f_k \rangle = 1, \qquad \forall k \in \Ind.
\end{equation}
Hence $\{h_k\}_{k \in \Ind}$ is a uniform tight frame.

We compute
%\begin{gather}
\begin{align}
\label{inf3}
\Big| |\langle f_k,f_l \rangle| - |\langle h_k,h_l \rangle| \Big|
& \le | \langle f_k,f_l \rangle - \langle h_k,h_l \rangle | \\
\label{inf4}
& \le | \langle f_k,f_l-h_l \rangle| + |\langle f_k-h_k,h_l \rangle | \\
\label{inf5}
& \le \|f_k\| \|f_l - h_l\| + \|h_l\| \|f_k- h_k\|,
\end{align}
%\end{gather}
where we have used the triangle inequality and the Cauchy-Schwarz
inequality. Note that
\begin{align}
\label{inf6}
\|f_l - h_l\| & = \|f_l - \sqrt{\rho} \SQI f_l \| \\
\label{inf7}
& \le \|(I - \sqrt{\rho} \SQI)\| \|f_l\| \\
\label{inf8}
& \le \max \Big\{|1-\sqrt{\frac{\rho}{A}}|, |1-\sqrt{\frac{\rho}{B}}|\Big\}.
\end{align}
Hence
\begin{equation}
\label{inf9}
\Big| |\langle f_k,f_l \rangle| - |\langle h_k,h_l \rangle| \Big|
\le 2 \max \Big\{|1-\sqrt{\frac{\rho}{A}}|, |1-\sqrt{\frac{\rho}{B}}|\Big\},
\end{equation}
and therefore
\begin{equation}
\label{inf10}
\underset{k,l \in \Ind; k \neq l}{\max} |\langle h_k,h_l \rangle|
\le \underset{k,l \in \Ind; k \neq l}{\max} |\langle f_k,f_l \rangle|
+ 2 \max \Big\{|1-\sqrt{\frac{\rho}{A}}|, |1-\sqrt{\frac{\rho}{B}}|\Big\}.
\end{equation}
\end{proof}

\remark 
(i) Although the canonical tight frame function $h_0$ do not have
to belong to $\Phi$, it is ``as close as possible'' to the function
$f_0 \in \Phi$. Indeed, under the assumptions of Theorem~\ref{th:tight} the 
(scaled) canonical tight frame $\tightframe$ generated by
$h_0 = \sqrt{\rho} \SQI f_0$ minimizes $\|f_0 - \phi_0\|$ among all tight 
frames $\{\phi_k\}_{k \in \Ind}$ in $\family(\Hsp,\Usp,L_2(\Hsp))$ 
(in fact among all possible tight frames), cf.~\cite{Han01} and for the 
case of Gabor frames~\cite{JS00}. However it is in general not true that 
$\tightframe$ also minimizes the maximal frame correlation 
$\max_{k,l}|\langle \phi_k,\phi_l \rangle|$ among all tight frames
$\{\phi_k\}_{k \in \Ind}$. \\
(ii) If the Grassmannian frame $\frame$ is already tight, then
the frame bounds satisfy $A=B=\rho$ and the second term in the
right-hand-side of~\eqref{tightgrass} vanished, as expected. \\
(iii) Frames that satisfy the assumptions of Theorem~\ref{th:tight} include
shift-invariant frames, Gabor frames and so-called geometrically
uniform frames (see~\cite{EB02} for the latter).

\subsection{An example: Grassmannian Gabor frames} \label{ss:gabex}

In this section we derive an example for Grassmannian frames in $\LtR$.
We consider Gabor frames in $\LtR$ generated by general lattices. 

Before we proceed we need some preparation. For $x, y \in \Rst$ we define 
the unitary operators of translation and modulation by $T_xf(t) = f(t-x)$, 
and $M_{\omega}f(t) = e^{2\pi i \omega t}f(t)$, respectively. Given a 
function $f \in \LtR$ we denote the time-frequency shifted function 
$f_{x,\omega}$ by
\begin{equation}
f_{x,\omega}(t) = e^{2\pi i\omega t}f(t-x).
\label{tfshift}
\end{equation}

A lattice $\Lambda$ of $\Rst^2$ is a discrete subgroup with compact
quotient. Any lattice is determined by its (non-unique) generator matrix 
$L \in GL(2,\Rst)$ via $\Lambda = L \Zst^2$. The volume of the lattice 
$\Lambda$ is $\vol(\Lambda) = \det(L)$.  

For a function (window) $g \in \LtR$ and a lattice $\Lambda$ in the
{\em time-frequency plane} $\Rst^2$ we define the corresponding Gabor system 
${\cal G}(g,\Lambda)$ by
\begin{equation}
{\cal G}(g,\Lambda) = \{M_{\omega}T_x g, \enspace  (x,\omega) \in \Lambda\}
\label{gabsys}
\end{equation}
Setting $\lambda = (x,\omega)$ we denote $\glam = M_{\omega}T_x g$.
If $\gab$ is a frame for $\LtR$ we call it a Gabor frame. As in remark~(ii) 
below Definition~\ref{def:red} we conclude that the redundancy of $\gab$ is 
$\rho=1/\vol(\Lambda)$. A necessary but by no means sufficient condition for 
$\gab$ to be a frame is $\vol(\Lambda) \le 1$, cf.~\cite{Gro01}.
It is clear that $\max_{\lambda \neq \lambda'}|\langle \glam,\glamp \rangle|$ 
will depend on the volume of the lattice, i.e., on the redundancy of the 
frame. The smaller $\vol(\Lambda)$ the larger 
$\max_{\lambda \neq \lambda'}|\langle \glam,\glamp \rangle|$.

One of the main purposes of Gabor frames is to analyze the
time-frequency behavior of functions~\cite{FS98}. To that end
one employs windows $g$ that are well-localized in time and frequency.
The Gaussian $\phi_{\sigma}(x)=(2/\sigma)^{\frac{1}{4}} e^{-\pi\sigma x^2},
\sigma >0$, is optimally localized in the sense that it minimizes
the Heisenberg Uncertainty Principle. Therefore Gabor frames using
Gaussian windows are of major importance in theory and applications.
Our goal is to construct Grassmannian Gabor frame generated by Gaussians.
Recall that $\gabg$ is a Gabor frame for $\LtR$ whenever 
$\vol(\Lambda) < 1$, see~\cite{Gro01}. 
Thus in the notation of Definition~\ref{def:grassframe} we consider
$\Hsp = \LtR$, $\Usp = \{T_x M_y, x,y \in \Lambda \,\, \text{with}\,\, 
\vol(\Lambda)=\rho\}$, and $\Phi=\{\phi_{\sigma} \, | \, \phi_{\sigma}(x) = 
(2\sigma)^{\frac{1}{4}} e^{-\pi x^2/\sigma}, \, \sigma > 0\}$.
That means for fixed redundancy $\rho$ we want to find
$\Lambdaopt$ among all lattices $\Lambda$ with
$\vol(\Lambda)=\rho$ and $\gaussopt$ among all Gaussians $\gauss$ such that 
\begin{equation}
\label{maxlatt}
\max_{\lambda\neq\lambda'}
|\langle (\gauss)_{\lambda},(\gauss)_{\lambda'} \rangle|
\end{equation}
is minimized.

Since $\hat{\phi}_{\sigma} = \phi_{1/\sigma}$ we can restrict our analysis 
to Gaussians with $\sigma=1$, as all other cases can be obtained by a 
proper dilation of the lattice. To simplify notation we write $\phi :=\phi_1$.

Since $T_x$ and $M_{\omega}$ are unitary operators there holds 
$|\langle \glam,\glamp \rangle| = |\langle g,g_{\lambda'-\lambda} \rangle|$ 
for any $g \in \LtR$. Furthermore $|\langle \phi,\gausslam \rangle|$ is 
monotonically decreasing with increasing $\|\lambda\|$ 
(where $\|\lambda\| = \sqrt{|x|^2+|\omega|^2}$)
due to the unimodality, symmetry, and Fourier invariance of $\phi$.
These observations imply that our problem reduces to
finding the lattice $\Lambdaopt$ of redundancy $\rho$ such that
$\max |\langle \phi,\gausslam \rangle|$ is minimized where
\begin{equation}
\lambda \in \big\{L e_1, L e_2, \enspace \text{with} \,\,
e_1 = [1,0]^T, e_2=[0,1]^T \big\}.
\label{lambda01}
\end{equation}
The ambiguity function of $f \in \LtR$ is defined as
\begin{equation}
Af(t,\omega) = \int \limits_{-\infty}^{+\infty} 
f(x+\frac{t}{2}) \overline{f(x-\frac{t}{2})} e^{-2\pi i \omega x} \, dx. 
\label{defamb}
\end{equation}
There holds $|\langle f,g \rangle|^2 = |\langle Af, Ag \rangle|$
for $f,g \in \LtR$. It follows from Proposition~4.76 in~\cite{Fol89} 
that $A\phi$ is rotation-invariant. Furthermore
$A \gausslam$ is rotation-invariant with respect to its ``center''
$\lambda = (x,\omega)$ which follows from~(4.7) and~(4.20) in
\cite{Gro01} and the rotation-invariance of $A\phi$.

Next we need a result from sphere packing theory. Recall that in
the classical sphere packing problem in $\Rst^d$ one tries to find the 
lattice $\Lambdaopt$ among all lattices $\Lambda$ in $\Rst^d$ that solves
\begin{equation}
\underset{\Lambda}{\max} 
\left\{\frac{\text{Volume of a sphere}}{\vol(\Lambda)}\right\}.
\label{sphere}
\end{equation}
For a given lattice $\Lambda$ the radius $r$ of such a sphere is 
\begin{equation}
r = \frac{1}{2}\big(\underset{
{\footnotesize \begin{matrix}
\lambda,\lambda' \in \Lambda \\
\lambda \neq \lambda' 
\end{matrix}}
}{\min} \{\|\lambda - \lambda'\|\}\big).
\label{radius}
\end{equation}
Hence solving~\eqref{sphere} is equivalent to solving
\begin{equation}
\underset{\Lambda}{\max}\big\{ \underset{
{\footnotesize \begin{matrix}
\lambda,\lambda' \in \Lambda \\
\lambda \neq \lambda 
\end{matrix}}
}{\min} \{\|\lambda - \lambda'\|\}\big\}
\qquad \text{subject to $\vol(\Lambda) = \rho$}
\label{sphere2}
\end{equation}
for some arbitrary, but fixed $\rho >0$. Obviously the minimum has to be 
taken only over adjacent lattice points.

Due to the rotation-invariance of $A\phi$ and $A\gausslam$ and since
$A\phi(x,\omega)$ is monotonically decreasing with increasing $(x,\omega)$
we see that the solution of
\begin{equation}
\underset{\Lambda}{\min}
\underset{\text{$\lambda$ as in~\eqref{lambda01}}}{\max}
\big\{ |\langle \phi,\phi_{\lambda} \rangle|\big\} \quad
\qquad \text{subject to $\vol(\Lambda) = \rho$}
\label{mininn}
\end{equation}
is identical to the solution of~\eqref{sphere2}. It is well-known
that the sphere packing problem~\eqref{sphere2} in $\Rst^2$ is
solved by the hexagonal lattice $\Lambda_{hex}$, see~\cite{CS93b}. Thus 
for given redundancy $\rho>1$ the Gabor frame ${\cal G}(\phi,\Lambda_{hex})$
is a Grassmannian frame, where the generator matrix of $\Lambda_{hex}$ is 
given by
\begin{equation}
L_{\text hex} =
\begin{bmatrix}
\frac{\sqrt{2}}{\sqrt[4]{3}\sqrt{\rho}} &
\frac{1}{\sqrt[4]{3}\sqrt{2\rho}} \\
0 & \frac{\sqrt[4]{3}}{\sqrt{2\rho}} 
\end{bmatrix}.
\label{hexlat}
\end{equation}

\remark (i) The result can be generalized to higher dimensions,
since the ambiguity function of a $d$-dimensional Gaussian is
also rotation-invariant. Hence a Grassmannian Gabor frame with
Gaussian window is always associated with the optimal lattice sphere packing
in $\Rst^{2d}$. However in higher dimensions explicit solutions
to the sphere packing problem are in general not known.~\cite{CS93b}. \\
(ii) If we define the Gaussian with complex exponent $\sigma = u+iv$
with $u >0$ (i.e., chirped Gaussians in engineering terminology) then it
is not hard to show that a properly chirped Gaussian associated with a
rectangular lattice yields also a Grassmannian Gabor frame.

\section{Applications of Grassmannian frames} \label{s:appl}

In this section we describe two applications of Grassmannian frames.
The first one concerns wireless communication and
involves infinite-dimensional Grassmannian frames, while the second
one concerns coding theory and involves mainly finite Grassmannian frames.
Another important application of finite Grassmannian frames is described 
in~\cite{HS02} where this concept is used to construct spreading sequences
for CDMA.

\subsection{Grassmannian Gabor frames and OFDM} \label{ss:OFDM}

We briefly describe an application of Grassmannian Gabor frames
in wireless communication. 
Orthogonal frequency division multiplexing (OFDM) has emerged
as attractive candidate for 4-th generation wireless communication
systems. We refer to~\cite{FAB95} for details about OFDM and
to~\cite{Str01} for details about the connection to Gabor theory. 
The transmission functions of an OFDM system are
\begin{equation}
\phi_{k,l}(t) = \phi(t - kT) e^{2\pi i lFt},   \quad k \in \Zst,
l = 0,\dots,M-1,
\label{ofdm}
\end{equation}
where $\phi$ is a given pulse shape. To minimize distortions caused
by additive white Gaussian noise the functions $\phi_{k,l}$ are usually
chosen to be mutually orthogonal. OFDM transmission
functions form a critically sampled or undersampled Gabor system, since 
perfect reconstruction at the receiver requires $TF \ge 1$. The duality
relations for Gabor systems provides the connection between Gabor frames
and undersampled Gabor systems~\cite{FS98}. 

The mobile wireless channel $\channel$ is time-dispersive as well as
frequency-dispersive which can lead to intersymbol interference (ISI) and
interchannel interference (ICI). When the pulses $\phi_{k,l}$ pass through
the channel $\channel$, the mutual orthogonality is lost.
Thus to mitigate ISI and ICI and to allow for a simple receiver
structure it is important that the pulse shapes $\phi_{k,l}$ are
well concentrated in time and frequency. Since the Balian-Low theorem
prohibits well-localized Gabor systems when $TF=1$ one usually chooses
$TF>1$. In this situation Gaussians are
a natural starting point for pulse shape design.

If the pulses are well-localized in time and frequency we have
\begin{equation}
\label{channel1}
|\langle \channel \phi_{k,l}, \channel \phi_{k \pm 1, l \pm 1} \rangle|
\gg |\langle \channel \phi_{k,l}, \channel \phi_{k \pm k',l \pm l'} \rangle|, 
\qquad \text{for}\,\,\, k',l' > 1,
\end{equation}
hence it is sufficient to consider the reduction of interference between
pulses that are adjacent in the time-frequency domain.

We can of course always increase the parameters $T$ and $F$ in order
to reduce interference, however this results in a lower data rate.
Therefore in~\cite{SB01} the OFDM scheme has been generalized by replacing the
rectangular time-frequency lattice $\{(kT,lF)\}_{k,l \in \Zst}$ by
a general lattice, yielding a Gabor system of the form~\eqref{gabsys}.
As follows from Subsection~\ref{ss:gabex} when using Gaussian pulses in
combination with a hexagonal-type lattice we can minimize
the interference between pulses by maximizing the distance
between adjacent lattice points without reducing the data rate. 
Gaussian pulses are not mutually orthogonal, but we can easily 
orthogonalize this system via the ``inverse square root method'', 
see~\cite{Str01}. The ''orthogonalized'' pulses will have a somewhat
larger interference as compared to Gaussian pulses, but this increase
in interference can be estimated conveniently via Theorem~\ref{th:tight}.
We refer to~\cite{SB01} for more details.

\subsection{Erasures, coding, and Grassmannian frames} \label{ss:coding}

\if 0
As mentioned in the introduction uniform tight frames possess
many properties that play an important role in applications. Tightness for
instance guarantees that additive noise that might distort frame
coefficients is not amplified in subsequent processing. Tight frames
also give rise to simple reconstruction formulas. When designing
filter banks, codes, or transmissions systems for wireless communication
it is often crucial that all ``building blocks'' (i.e, filters, code words, or 
transmission pulses) have the same energy. In other words each 
frame element has the same norm. This property together with tightness
singles out uniform tight frames as particularly important class of frames.
\fi

Recently finite frames have been proposed for multiple description coding for 
erasure channels, see~\cite{GKK01,GK01,CK02}. We consider the following setup.
Let $\finframe$ be a frame in $\Est^m$. As in~\eqref{frameanalysis}
and~\eqref{framesynthesis} we denote the associated analysis and synthesis 
operator by $T$ and $T^{\ast}$ respectively. Let $f \in \Est^m$ represent
the data to be transmitted. We compute $y=Tf \in \Est^N$ and send $y$ over the
erasure channel. We denote the index set that corresponds to the erased 
coefficients by $\erind$ and the surviving coefficients are indexed by the set
$\remind$. Furthermore we define the $N \times N$ {\em erasure matrix} $Q$ via
\begin{equation}
\label{Qmat}
Q_{kl} = 
\begin{cases}
0 & \text{if $k \neq l$,} \\
0 & \text{if $k = l$ and $k \in {\erind}$,} \\
1 & \text{if $k = l$ and $k \in {\remind}$.} \\
\end{cases}
\end{equation}
Let $\eps$ represent additive white Gaussian noise (AWGN) with zero mean and
power spectral density $\sigma^2$. The data vector arriving at the receiver 
can be written as $\yerr := Qy +\eps$.

The frame $\finframe$ is robust against $e$ erasures, 
if $\erframe$ is still a frame for $\Est^m$ for any index set 
$\remind \subset \{0,\dots,N-1\}$ with $|\remind| \ge N-e$. In this
case standard linear algebra implies that $f$ can be exactly reconstructed 
from $\yerr$ in the absence of noise\footnote{For instance so-called 
harmonic frames 
are robust against up to $N-m$ erasures~\cite{GKK01}, which does not come as 
a surprise to those researchers who are familiar with Reed-Solomon codes or 
with the fundamental theorem of algebra.}.

In general, when we employ a minimum mean squared (MMSE) receiver we compute 
the (soft) estimate
\begin{equation}
\frec = (\TES \TE+\sigma^2 I_m)^{-1} \TES \yerr, 
\label{errrec}
\end{equation}
where $\TE$ is the analysis operator of the frame $\erframe$.
This involves the inversion of a possibly large matrix (no matter 
if noise is present or not) that can differ from one transmission to the
next one. The costs of an MMSE receiver may be prohibitive in time-critical 
applications. Therefore one often resorts to a matched filter receiver
which computes the estimate
\begin{equation}
\frec = c T^{\ast} \yerr,
\label{matchedrec}
\end{equation}
where $T$ is the analysis operator of the original frame $\finframe$
and $c>0$ is a scaling constant.
The advantage of an MMSE receiver is the better error performance while
a matched filter receiver can be implemented at lower computational cost.

Robustness against the maximal number of erasures is not the only 
performance criterion when designing frames for coding. Since any 
transmission channels is subject to AWGN, it is important that the noise 
does not get amplified during the transmission process. Yet another 
criterion is ease of implementation of the receiver. It is therefore
natural to assume $\finframe$ to be a uniform tight frame,
since in case of no erasures (i)~the MMSE receiver coincides with the matched 
filter receiver and (ii)~AWGN does not get amplified during transmission.

Our goal in this section is to design a uniform tight frame such that the 
performance of a matched filter receiver is maximized in presence of
erasure channels. In other words the  approximation error 
$\|f - \tilde{f}\|$ is minimized, where $\tilde{f}$ is computed via
a matched filter receiver, i.e., $\tilde{f} = \frac{m}{N}\TT \tilde{y}$
with $\tilde{y}= Qy +\eps$.

We estimate the reconstruction error via
\begin{align}
\|f -\tf\| = & \|f - \frac{m}{N}\TT (Q T f + \eps)\| \\
\le & \frac{m}{N}\|\TT T f- \TT Q T f\| + \frac{m}{N}\|\TT \eps\| \\
\le & \frac{m}{N} \|\TT P T\|_2 \| f\| + \sigma, 
\label{est1}
\end{align}
where we used the notation $P=I-Q$. Since $P$ is an orthogonal projection 
and therefore satisfies $P P^{\ast}=P$ there holds
\begin{gather}
\| \TT P T\|_2 = \|(PT)^{\ast} PT\|_2 = \|PT (PT)^{\ast}\|_2
= \| P T\TT  P\|_2.
\label{est2}
\end{gather}
Hence we should design our tight uniform frame $\frame$ such 
that $\| P T\TT  P\|_2$ is minimized, where the minimum
is taken over all matrices $P=I-Q$, with $Q$ as defined in~\eqref{Qmat}.

Recall that $T \TT =\{\langle f_l, f_k \rangle\}_{k,l \in \Ind}$,
hence $P T\TT P = \{\langle f_l, f_k \rangle\}_{k,l \in \Indo}$.
Furthermore
\begin{equation}
\|P T\TT P\|_2 \le \sqrt{\|P T\TT P\|_{\infty} \|P T\TT P\|_1} = 
\max_{k \in \Indo}\sum_{l \in \Indo}|\langle f_k, f_l \rangle|.
\label{est3}
\end{equation}
This suggests to look for frames for which 
$\underset{k,l, k\neq l}{\max}|\langle f_k, f_l \rangle|$
is minimized, in other words to look for optimal Grassmannian frames.

\remark (i) In case of one erasure it has been shown in~\cite{GKK01}
(in case of unknown $\sigma$) that uniform tight frames are optimal with 
respect to minimizing the influence of AWGN when using the MMSE receiver. 
In case of one erasure uniform tight frames also minimize
the reconstruction error when using a matched filter receiver. \\
(ii) Holmes and Paulsen have shown that Grassmannian frames are optimal 
with respect to up to two erasures~\cite{HP02}. This can be easily seen by
minimizing the operator norm of the matrix $P T T^{\ast}P$, which
in this case reduces {\em exactly} to the problem of 
minimizing $\max |\langle f_k,f_l \rangle|$ for all $k,l$ with $k \neq l$. \\
(iii) There is strong numerical evidence that the optimal Grassmannian frames
of part~(b) in Corollary~\ref{confgrass} are even optimal for three erasures 
(however this is not the case for the frames constructed in part~(a)).\\
(iv) Grassmannian frames are in general not robust against $N-m$ erasures
if $N>m+2$.

\bigskip

\noindent
{\bf Example:} We elaborate further an example given in~\cite{GKK01}, where
the authors consider the design of multiple description coding frames 
$\finframe$ in $\Est^m$ with $m=3$ and $N=7$. As in Examples~4.2 and~4.3 
in~\cite{GKK01} we consider an erasure channel with AWGN, but unknown
noise level. Without knowledge of $\sigma$ the 
the reconstruction formula of the MMSE receiver simplifies to
$\frec = (\TES \TE)^{-1} \TES \yerr$. Standard numerical analysis
tells us that the smaller the condition number of $\TES \TE$ the smaller
the amplification of the noise in the reconstruction \cite{GL96}. We
therefore compare the condition number of different uniform tight frames
for $m=3, N=7$ after up to four frame elements have been randomly removed.

We consider three types of uniform tight frames. The first frame is an optimal
complex-valued Grassmannian frame, the second one is constructed by taking
the first three rows of a $7 \times 7$ DFT matrix and using the columns
of the resulting (normalized) $3 \times 7$ matrix as frame elements
(this is also called a harmonic frame in~\cite{GKK01}).
The last frame is a randomly generated uniform tight frame.
Since all three frames are uniform tight, they show identical performance
for one random erasure, the condition number of the frame operator
in this case is constant 1.322. 
Since the frames are of small size, we can easily compute the condition
number for all possible combinations of two, three, and four erasures. 
We then calculate the maximal and mean average condition number for
each frame. As can be seen from the results in Table~\ref{table1} 
the optimal Grassmannian frame outperforms
the other two frames in all cases, except for the average condition 
number for two erasures, where its condition number is slightly larger.
This example demonstrates the potential of Grassmannian frames
for multiple description coding.

\begin{table}[h]
\label{table1}
\begin{center}
\begin{tabular}{|c|c|c|c|c|c|c|}
\hline
   & \multicolumn{2}{c|}{2 erasures}
& \multicolumn{2}{c|}{3 erasures}   & \multicolumn{2}{c|}{4 erasures}\\
\hline
Cond.number & mean & max & mean & max & mean & max \\
\hline
Grassmannian frame & 1.645& 1.645 & 2.044& 2.188 & 3.056& 3.630 \\
DFT-submatrix frame   & 1.634& 1.998 & 2.198& 3.602 & 4.020& 8.580 \\
Random uniform tight & 1.638& 1.861 & 2.095& 3.792 & 3.570&12.710 \\
\hline 
\end{tabular}
\label{tablecomp}
\caption{Comparison of mean average and maximal condition number of
frame operator in case of two, three, and four random erasures. We compare
an optimal Grassmannian frame, a DFT-based uniform tight frame, and
a random uniform tight frame. The Grassmannian frame shows the best overall
performance.} 
\end{center}
\end{table}

%\bibliographystyle{plain}
%\bibliography{nuhag,linalg,gabor,mathbook,tf-zentral,code,thomas}

\end{document}